\documentclass[10pt]{article}

\usepackage{amscd,amsmath, amssymb, fancyhdr, mathabx}
\usepackage[backref=page]{hyperref}

\renewcommand*{\backrefalt}[4]{%
	\ifcase #1 (Not cited.)%
	\or        (Cited on page~#2.)%
	\else      (Cited on pages~#2.)%
	\fi}

\hypersetup{
	colorlinks   = true,
	citecolor    = magenta,
	linkcolor    =blue,
	urlcolor     =magenta	
}

\numberwithin{equation}{section}

\newcommand{\version}{version 1.0,\ \ May 26, 2026}

\def\eqref#1{(\ref{#1})}

\def\1{\sqrt{-1}\:}

\newcommand{\cntrct}                
{\hspace{2pt}\raisebox{1pt}{\text{$\lrcorner$}}\hspace{2pt}}

\renewcommand{\tilde}{\widetilde}
\renewcommand{\bar}{\overline}
\renewcommand{\phi}{\varphi}
\renewcommand{\epsilon}{\varepsilon}
\renewcommand{\geq}{\geqslant}

\newcommand{\Bim}{\mathrm{Bir}}
\newcommand{\Aut}{\mathrm{Aut}}
\newcommand{\MBM}{\mathrm{MBM}}
\newcommand{\NS}{\mathrm{NS}}
\newcommand{\Mon}{\mathrm{Mon}}

\newcommand{\Hom}{\operatorname{Hom}}

\newcommand{\Diff}{\operatorname{Diff}}

\newcommand{\rk}{\operatorname{rk}}

\newcommand{\bbL}{\mathbb{L}}
\newcommand{\bbM}{\mathbb{M}}
\newcommand{\bbI}{\mathbb{I}}

\newcommand{\bbN}{\mathbb{N}}
\newcommand{\bbZ}{\mathbb{Z}}
\newcommand{\bbQ}{\mathbb{Q}}
\newcommand{\bbR}{\mathbb{R}}
\newcommand{\bbC}{\mathbb{C}}

\newcommand{\bbP}{\mathbb{P}}

\newcommand{\CC}{\mathcal{C}}
\newcommand{\DD}{\mathcal{D}}

\newcommand{\MM}{\mathcal{M}}
\newcommand{\NN}{\mathcal{N}}

\newcommand{\UU}{\mathcal{U}}
\newcommand{\VV}{\mathcal{V}}

\newcommand{\XX}{\mathcal{X}}
\newcommand{\YY}{\mathcal{Y}}

\newcommand{\SO}{\mathrm{SO}}

\newcommand{\rmO}{\mathrm{O}}
\newcommand{\rmW}{\mathrm{W}}
\newcommand{\rmR}{\mathrm{R}}
\newcommand{\st}{\enskip |\enskip}

\newcounter{Mycounter}[section]
\newcounter{lemma}[section]
\renewcommand{\thelemma}{{Lemma \thesection.\arabic{lemma}}}
\newcommand{\lemma}{%
    \setcounter{lemma}{\value{Mycounter}}
    \refstepcounter{lemma}
    \stepcounter{Mycounter}
    {\noindent \bf \thelemma:\ }}

\newcounter{claim}[section]
\newcounter{sublemma}[section]
\newcounter{corollary}[section]
\renewcommand{\thecorollary}{{Corollary \thesection.\arabic{corollary}}}
\newcommand{\corollary}{%
    \setcounter{corollary}{\value{Mycounter}}
    \refstepcounter{corollary}
    \stepcounter{Mycounter}
    {\noindent \bf \thecorollary:\ }}

\newcounter{theorem}[section]
\renewcommand{\thetheorem}{{Theorem \thesection.\arabic{theorem}}}
\newcommand{\theorem}{%
    \setcounter{theorem}{\value{Mycounter}}
    \refstepcounter{theorem}
    \stepcounter{Mycounter}
    {\noindent \bf \thetheorem:\ }}

\newcounter{conjecture}[section]
\newcounter{proposition}[section]
\renewcommand{\theproposition}
      {{Proposition \thesection.\arabic{proposition}}}
\newcommand{\proposition}{%
    \setcounter{proposition}{\value{Mycounter}}
    \refstepcounter{proposition}
    \stepcounter{Mycounter}
    {\noindent \bf \theproposition:\ }}

\newcounter{definition}[section]
\renewcommand{\thedefinition}
      {{Definition~\thesection.\arabic{definition}}}
\newcommand{\definition}{%
    \setcounter{definition}{\value{Mycounter}}
    \refstepcounter{definition}
    \stepcounter{Mycounter}
    {\noindent \bf \thedefinition:\ }}

\newcounter{example}[section]
\newcounter{remark}[section]
\renewcommand{\theremark}{{Remark \thesection.\arabic{remark}}}
\newcommand{\remark}{%
    \setcounter{remark}{\value{Mycounter}}
    \refstepcounter{remark}
    \stepcounter{Mycounter}
    {\noindent \bf \theremark:\ }}

\newcounter{problem}[section]
\newcounter{question}[section]
\newcommand{\proof}{\noindent{\bf Proof:\ }}

\makeatletter

\@addtoreset{equation}{section} \@addtoreset{footnote}{section}
\makeatother

\def\blacksquare{\hbox{\vrule width 5pt height 5pt depth 0pt}}
\def\endproof{\blacksquare}

\begin{document}
\begin{center}
{\LARGE\bf
Reflective lattices and hyperk\"ahler manifolds
\\[4mm]
}

	{Ekaterina Amerik\footnote{Partially supported by the HSE University Basic Research Program},
Andrey Soldatenkov\footnote{Partially supported by 
FAPESP grant 2024/23819-0}, 
Misha Verbitsky\footnote{Partially supported by 
FAPERJ grant SEI-260003/000410/2023 and CNPq - Process 310952/2021-2. \\


 }}

\end{center}

{\small \hspace{0.10\linewidth}
\begin{minipage}[t]{0.85\linewidth}
{\bf Abstract:}
Using the results of Nikulin and Vinberg on the
groups of isometries generated by reflections, we construct a subvariety
called the Nikulin--Vinberg locus in the moduli space of
polarized hyperk\"ahler manifolds. It is obtained as a finite union
of components of higher Noether--Lefschetz loci which parameterize manifolds with 
certain special N\'eron--Severi lattices. The Nikulin--Vinberg locus
is the closure of the set of hyperk\"ahler manifolds with Picard number $\geq 3$
which have finite groups of birational automorphisms. 
Using this construction and a refinement of an argument by Oguiso, we show that any non-trivial family
of projective deformations of a hyperk\"ahler manifold with $b_2(M)\geq 6$
has a dense set of fibers which have an infinite group of
birational automorphisms.
\end{minipage}
}

\tableofcontents

\section{Introduction}

The goal of this paper is to construct a certain family of special
subvarieties in the moduli spaces of polarized hyperk\"ahler manifolds,
and to discuss the relation of those subvarieties with the properties
of birational automorphism groups of the corresponding hyperk\"ahler
manifolds. Let us state our main result. We will denote by $M$ a hyperk\"ahler
manifold in a fixed deformation class, by $h\in H^2(M, \bbZ)$ the class
of a polarization and by $\MM_h$ the moduli space of hyperk\"ahler
manifolds polarized by $h$. We denote by $[M]\in \MM_h$ the point 
corresponding to $M$,
by $\Aut(M)$ and $\Bim(M)$ the groups of biregular and birational automorphisms
of $M$, respectively.
We will recall all necessary definitions below in Section \ref{sec_hk}.

\hfill

\theorem\label{thm_main} Assume that $b_2(M)\ge 6$. There exists a finite family of
special subvarieties $\XX_i\subset \MM_h$, $i\in \{1,\ldots,m\}$
and a countable family of special divisors $\YY_j\subset \MM_h$, $j\in \bbN$
that have the following properties.
\begin{enumerate}
\item Assume that $M$ is a hyperk\"ahler manifold with $[M]\in \MM_h$,
$\rk\NS(M)\ge 3$ and $|\Bim(M)| < \infty$. Then there exists a unique
$i\in \{1, \ldots, m\}$ such that:
\begin{enumerate}
	\item $[M]\in \XX_i$;
	\item For any $M'$ such that $[M']$ is a very general point of $\XX_i$,
	there exists an isometry between $\NS(M)$ and $\NS(M')$.
\end{enumerate}
\item The family of divisors $\{\YY_j\}_{j\in \bbN}$ is strongly dense
in $\MM_h$ in the sense of \ref{def_density}.
\item For every $j\in \bbN$ there exists a Zariski-open non-empty
subset $\YY_j^\circ\subset \YY_j$ with the following property.
Assume that $M$ is a hyperk\"ahler manifold such that $[M]\in \YY_j^\circ$
for some $j\in \bbN$. Then $\Aut(M)$ contains an element of infinite order.
\end{enumerate}

\hfill

The first part of the theorem above is based on a deep result of Nikulin
\cite{_Nikulin_} and Vinberg \cite{_Vinberg81_}
on the finiteness of the number of reflective Lorentzian
lattices of rank at least three. That result was one of the major achievements
within a vast research program initiated by Vinberg \cite{_Vinberg67_},
dedicated to the study of discrete groups of isometries generated by
reflections in hyperbolic spaces.
We will therefore call the union $\cup_i \XX_i\subset \MM_h$ of the subvarieties
appearing in the above theorem the Nikulin--Vinberg locus. We refer to \cite{_Belolip_}
for an overview of more recent developments related to Vinberg's program.

 The first part of
the theorem can be briefly (and imprecisely) summarized by saying that the
polarized hyperk\"ahler manifolds with Picard rank at least three having
only finitely many birational automorphisms lie in the Nikulin--Vinberg locus,
which is a proper closed subvariety of the moduli space.

The second and third parts of the theorem are 
based on an explicit (and elementary) construction
of binary anisotropic\footnote{A quadratic lattice is
  called anisotropic if it does not contain a non-zero vector $x$ with
  $q(x)=0$.}
 Lorentzian lattices of arbitrarily 
large discriminant, the construction that we will explain below. 
The above theorem has the following corollary.

\hfill

\corollary\label{cor_main}
Let $\pi\colon \MM \to B$ be a smooth non-isotrivial
family of polarized hyperk\"ahler manifolds with $b_2\ge 6$
over a complex base $B$.
Then there exists a dense subset $B'\subset B$ such that for $t\in B'$
the group $\Bim(M_t)$ is infinite, where $M_t = \pi^{-1}(t)$.

\hfill

This corollary was first proven for K3 surfaces by Oguiso in \cite{_Oguiso_},
and then extended to all known deformation types of hyperk\"ahler manifolds in \cite{_Denisi_}.
The result in fact holds in greater generality
and does not rely on any information specific to the known families of hyperk\"ahler manifolds.
Both the above theorem and its corollary readily extend, for example, to hyperk\"ahler orbifolds.
Our proof is purely lattice-theoretic. The only technical condition that we impose concerns
$b_2$, but we believe that it should be possible to get rid of this condition.
We also note that the first part of the theorem does hold without assuming $b_2\ge 6$.

\hfill

\remark\label{aut??}
So far we don't know whether an analogue of \ref{cor_main} holds with $\Bim(M_t)$ replaced by $\Aut(M_t)$.
One way to prove such stronger form of \ref{cor_main} is by showing that the preimage in $B$
of the locus $\cup_j \YY_j^\circ$ under the period map is dense in $B$. Here $\YY_j^\circ$ are the open
subsets of the divisors $\YY_j$ appearing in part 3 of \ref{thm_main}. General considerations
along the lines of the Zilber--Pink conjecture indicate that this could be a reasonable expectation,
at least in the case when the base $B$ is algebraic.

\hfill

\subsection{Lattices}\label{sec_lat}

An integral lattice is a free $\bbZ$-module $\Lambda\simeq \bbZ^r$ together with a non-degenerate
symmetric bilinear form $(\cdot,\cdot)\colon \Lambda\otimes\Lambda \to \bbZ$.
The corresponding quadratic form will be denoted by $q$, so that $q(v) = (v,v)$
for any $v\in\Lambda$.
We will also extend $(\cdot,\cdot)$ to a $\bbQ$-bilinear form on the vector
space $\Lambda\otimes \bbQ$.
The dual lattice $\Lambda^* = \Hom(\Lambda, \bbZ)$ is embedded into $\Lambda\otimes\bbQ$ via the form
$(\cdot,\cdot)$, so that
$$
\Lambda^*= \{v\in \Lambda\otimes \bbQ\st \forall u\in\Lambda \,\,\, (u,v)\in \bbZ\},
$$
and the quotient $\Lambda^*/\Lambda$ is called the discriminant group, denoted $d(\Lambda)$.
The exponent of the lattice is the smallest positive
integer  $e(\Lambda)$ which annihilates the discriminant group.

A vector $v\in \Lambda\otimes \bbQ$ with $q(v)\neq 0$ 
(such $v$ is called non-isotropic)
defines a reflection
$s_v$, that is the following involution of $\Lambda\otimes\bbQ$:
$$
s_v(u) = u - 2\frac{(u, v)}{(v,v)} v.
$$
A primitive non-isotropic element $v\in \Lambda$ is called a root if $s_v(\Lambda) = \Lambda$. Equivalently,
$v$ is a root if $q(v)$ divides $2 (u, v)$ for all $u\in \Lambda$. It is 
well-known (and follows directly from the definitions) that in this case $q(v)$
divides $2 e(\Lambda)$, see e.g. \cite[section 1]{_Scharlau_}. Any vector $v\in \Lambda$
with $0 < |q(v)| \le 2$ is automatically a root. Let us denote by $\rmR(\Lambda)$ the set
of all roots $v\in \Lambda$ with $q(v) < 0$. Note our convention for the sign of the roots:
the condition $q(v) < 0$ is not standard and not necessary for studying lattices in general, but
for our applications only the negative roots will be relevant.

Given a lattice $\Lambda$, consider the group $\rmO(\Lambda)$ of orthogonal transformations
of $\Lambda\otimes\bbQ$ preserving $\Lambda$. The reflections $s_v$, where $v\in \rmR(\Lambda)$, generate a subgroup $\rmW(\Lambda)\subset \rmO(\Lambda)$ called
the Weyl group of $\Lambda$. If the lattice $\Lambda$ has no negative roots, the Weyl group is trivial.

A lattice $\Lambda$ is called Lorentzian if the quadratic form $q$
has signature $(1, r - 1)$. Note that in the literature dedicated
to the theory of lattices another convention is often used: the Lorentzian signature is $(r-1, 1)$.
Our convention is more convenient in the context of algebraic geometry.

A lattice $\Lambda$ is called reflective if $W(\Lambda)$
is of finite index in $O(\Lambda)$. We refer to \cite{_Belolip_}
for an overview of the results related to the classification of reflective Lorentzian lattices.

\subsection{Moduli spaces of polarized hyperk\"ahler manifolds}\label{sec_hk}

By a hyperk\"ahler (also known as irreducible holomorphically symplectic
or IHS) manifold $M$ we will mean a compact simply connected complex $2n$-dimensional
manifold that admits a Riemannian metric of holonomy $\mathrm{Sp}(n)$,
see \cite{_Beauville_}, \cite{_Markman:Survey_} and \cite{V1} for the
definitions and the foundations of the theory of such manifolds.

Let $\Lambda = H^2(M, \bbZ)$ be the second cohomology of a hyperk\"ahler
manifold $M$ equipped with the Beauville--Bogomolov--Fujiki (BBF) form $q$.
The latter has signature $(3, b_2(M) - 3)$, and we will assume that it is
rescaled to be integral and primitive. All deformations of $M$ that admit
K\"ahler metrics will also be hyperk\"ahler manifolds with the same $\Lambda$
and $q$. 

From now on we fix a deformation class of hyperk\"ahler manifolds and
denote by $M$ one of the manifolds in this deformation class.
The moduli theory of the deformations of $M$ is governed by the period
domains attached to $\Lambda$. Let us recall their definitions.

Given a primitive sublattice $\bbL\subset \Lambda$ of signature $(1, k)$,
its orthogonal complement $\bbL^\perp$ in $\Lambda$ has signature $(2, b_2(M) - k)$.
Let
$$
\DD_\bbL = \{ x \in \bbP(\bbL^\perp \otimes \bbC)\st q(x) = 0,\, q(x, \bar{x}) > 0\}.
$$
This is a Hermitian symmetric domain whose points parametrize hyperk\"ahler
manifolds $M$ with $\bbL\subset \NS(M)$, where the latter is the N\'eron--Severi
group of $M$.
For a very general (that is, lying outside of a countable union
of proper subvarieties) point of $\DD_\bbL$ we have equality $\bbL = \NS(M)$.
We will call $\DD_\bbL$ the period domain of $\bbL$-polarized hyperk\"ahler manifolds.
Note that the condition on the signature of $\bbL$ guarantees that all $\bbL$-polarized
hyperk\"ahler manifolds are projective. In particular, if $\bbL = \bbZ h$ for some
primitive $h\in \Lambda$ with $q(h)> 0$, we will write $\DD_h$ instead of $\DD_\bbL$ and call it
the period domain of $h$-polarized hyperk\"ahler manifolds.
We will denote the orthogonal complement $h^\bot$ in $\Lambda$
by  $\Lambda_h$.

Clearly, given two primitive sublattices $\bbL_1 \subset \bbL_2$, we have
$\DD_{\bbL_2} \subset \DD_{\bbL_1}$. In particular, if $h\in \bbL$, we
have $\DD_{\bbL} \subset \DD_h$. We call such $\DD_\bbL$ a special subvariety
of the $h$-polarized period domain $\DD_h$.

The mapping class group $\Gamma=\Diff(M)/\Diff_0(M)$
acts on the Teichm\"uller space parametrizing the complex structures
on $M$ , see \cite{V1}. Let $\Gamma^M$ denote 
its subgroup preserving the connected component of the
Teichm\"uller space containing $M$.
Recall that the image of $\Gamma^M$  in $\rmO(\Lambda)$
under the natural action on the second cohomology is a finite index subgroup
of $\rmO(\Lambda)$ called the monodromy group, denoted $\Mon(M)$. For $\bbL\subset \Lambda$
we will denote by $\Mon_\bbL(M)\subset \Mon(M)$ the (pointwise) stabilizer
of the sublattice $\bbL$.

The image of $\Mon_\bbL(M)$ is of finite index in $\rmO(\bbL^\perp)$.
It is well-known (see \cite{Bai-Bo}) that the quotient $\MM_\bbL = \DD_{\bbL} / \Mon_\bbL(M)$
is a quasi-projective variety. We will call this quotient the moduli space of
$\bbL$-polarized hyperk\"ahler manifolds. When $\bbL=\bbZ h$, the moduli space of
$h$-polarized hyperk\"ahler manifolds will be denoted by
$\MM_h =\DD_h / \Mon_h(M)$. It comes with
the collection of special subvarieties $\MM_\bbL$ corresponding to the sublattices
$\bbL$ containing $h$.

Recall how one can describe the (bimeromorphic) automorphism group of $M$
in terms of $\NS(M)$. There exists a $\Mon(M)$-invariant collection of second
cohomology classes in $\Lambda$ called MBM classes, see \cite{_AV:MK_}. We will denote them
by $\MBM\subset \Lambda$. Recall from \cite{_AV:Orbits_} that $\MBM$ consists of a
finite number of $\Mon(M)$-orbits, in particular the BBF squares of the MBM classes are bounded.

Let $\bbL = \NS(M)$ and $\MBM_\bbL = \MBM \cap \bbL$. Denote by $\CC\subset \bbL\otimes \bbR$
the connected component of the positive cone that contains the K\"ahler classes. Then
for $x\in \MBM_\bbL$ the hyperplanes $x^\perp$ cut $\CC$ into a locally finite collection
of open subcones called K\"ahler chambers. The automorphism group of $M$ is isomorphic to
the subgroup of $\Mon(M)$ that preserves $\bbL$, the K\"ahler chamber of $M$ and the period of $M$.

The divisorial MBM classes on $M$ are elements of $\MBM_\bbL$ which are stably prime exceptional (\cite{_Markman:Survey_}, section 6; see also \cite{AV-RC}). Let $W^+(\bbL)$ be the subgroup of $W(\Lambda)$ that is
generated by reflections in divisorial MBM classes. 
Recall from \cite{_Markman:Survey_}  that we
have $W^+(\bbL) \subset \Mon(M)$. Since the roots that we are considering are always $q$-negative,
the action of $W^+(\bbL)$ on $\bbL\otimes \bbR$ preserves the cone $\CC$. The fundamental
domain of $W^+(\bbL)$-action on $\CC$ that contains a K\"ahler class of $M$ is called the
bimeromorphic K\"ahler chamber of $M$. The group of bimeromorphic automorphisms of $M$ is
isomorphic to the subgroup of $\Mon(M)$ that preserves the bimeromorphic K\"ahler chamber and the 
Hodge structure of $M$ (in particular,
$\bbL=\NS(M)$ is preserved automatically). It is finite exactly when $W^+(\bbL)$ is of finite index in the orthogonal group (\cite{_Markman:Survey_}, theorem 6.18; see also \cite{Vinberg-units} for general lattice-theoretic 
treatment).

\section{Subvarieties in the moduli space}\label{sec_subvar}

\subsection{The Nikulin--Vinberg locus}\label{sec_nikvin}

Here we will construct a family of subvarieties $\XX_i\subset \DD_h$
appearing in the first part of \ref{thm_main}. We fix a positive class
$h\in \Lambda$ with $q(h)=d$ that defines a polarization.

We will use the following notation: given a lattice $\bbL$ we denote by $\bbL(k)$
the same lattice with the form $(\cdot,\cdot)$ rescaled by an integer $k\neq 0$.
Also, let $\bbI$ denote the rank one lattice with $q(v) = v^2$ for $v\in \bbZ$.

\hfill

\proposition\label{prop_NV}
There exists a finite collection of lattices $\bbM_i$, $i\in \{1,\ldots, m_0\}$
with the following property. Suppose that $M$ is an $h$-polarized hyperk\"ahler
manifold such that $\rk \NS(M) \ge 3$ and $\Bim(M)$ is finite. Then $\NS(M)\cap \Lambda_h\simeq \bbM_i$
for some $i$.

\hfill

\proof
We denote by $\NN\VV$ the set of isomorphism classes of integral reflective Lorentzian 
lattices of rank at least three that are unscaled (i.e. such that the Gram matrix
of the lattice is indivisible). It follows form \cite[Theorem 5.2.1]{_Nikulin_} and
\cite[Theorems 2 and 4]{_Vinberg81_} that the set $\NN\VV$ is finite.

Recall from Section \ref{sec_hk} that under our assumptions on $M$ the orthogonal group
of the lattice $\NS(M)$ contains a finite index subgroup generated by reflections
in the roots that are also MBM classes. It follows that the lattice $\NS(M)$ is reflective,
though possibly not unscaled. Recall also that the squares of the MBM classes are
negative and bounded, the bound being deformation-invariant. Let $-N$ be the lower bound
for the squares of the MBM classes. Define
$\NN\VV^+ = \cup_{k = 1}^{N} \{\bbL(k)\st \bbL\in \NN\VV\}$.
The latter is a finite set of isomorphism classes of lattices, and we
see that $\NS(M)\in \NN\VV^+$.

To conclude, recall e.g. from \cite[Satz 30.2]{_Kneser:Quadratische_}
that for any lattice $\bbL$ there exist only finitely many primitive embeddings
$\bbI(d)\hookrightarrow \bbL$, up to the action of $\rmO(\bbL)$.
For every $\bbL \in \NN\VV^+$ and for every such primitive embedding
$\bbI(d)\hookrightarrow\bbL$ consider its orthogonal complement
$\bbM = \bbI(d)^\perp\subset \bbL$. Then there are, up to isomorphism,
only finitely many lattices that can occur as such orthogonal complements.
We denote these lattices $\bbM_i$, $i\in \{1,\ldots, m_0\}$.
Since $\NS(M)\in \NN\VV^+$ and $M$ is $h$-polarized, with $q(h) = d$,
we see that $\NS(M)\cap \Lambda_h$ is one of the orthogonal complements
constructed above.
\endproof

\hfill

We now define the family of subvarieties that form the Nikulin--Vinberg
locus in the moduli space. Recall from Section \ref{sec_hk} that $\Mon_h\subset \rmO(\Lambda_h)$
is a finite index subgroup. Combining this fact with \cite[Satz 30.2]{_Kneser:Quadratische_},
we conclude that for each of the lattices $\bbM_i$, $i\in \{1,\ldots, m_0\}$ from \ref{prop_NV}
there exist finitely many (up to the action of $\Mon_h$) primitive embeddings
$\bbM_i\hookrightarrow \Lambda_h$. Let us denote by $\bbM'_j\subset \Lambda_h$, $j\in \{1,\ldots, m\}$ the images of such embeddings that represent every of the $\Mon_h$-orbits,
for all $\bbM_i$ from \ref{prop_NV}.

\hfill

\definition\label{def_NV}
Let $\XX_i\subset \MM_h$, $i\in \{1,\ldots, m\}$ be the
image of $\DD_{\bbM'_i}$ under the quotient map $\DD_h\to \MM_h$.
The union of all $\XX_i$ is called the Nikulin-Vinberg
locus of $M$.


\subsection{Binary Lorentzian lattices}\label{sec_binlat}

Given $n\in \bbZ$, we will say that a lattice $\bbL$ represents $n$ if there exists
a non-zero $v\in \bbL$
with $q(v) = n$. For a lattice $\bbL$ let us denote by $\mu(\bbL)$ the following number
$$
\mu(\bbL) = \mathrm{max}\big( \{ q(v) \st v\in \bbL \}\cap \bbZ_{<0} \big),
$$
i.e. $\mu(\bbL)$ is the largest negative integer represented by $\bbL$.

Given two lattices, the symbol $\oplus$ will be used to denote their orthogonal direct sum.
A binary lattice is a lattice of rank two. We will consider
binary Lorentzian lattices $\bbL_a=(\bbZ^2, q)$, where
 $q=\begin{pmatrix} 1&0\\ 0& -a\end{pmatrix}$ for some integer $a > 0$.

\hfill

\proposition\label{prop_bin_lattices}
Let $n > 0$ be a fixed integer. Then for any integer $b>0$ there exist infinitely
many integers $a > 0$ such that $\bbL_{ab}$ does not contain isotropic vectors and
$\mu(\bbL_{ab}) < -n$, i.e. such that the lattice $\bbL_{ab}$ does
not represent $0, -1, -2, \ldots, -n$.

\hfill

\proof
For any $k\in \bbZ$ the condition that $\bbL_{ab}$ represents $-k$ means that
there exist $x, y\in \bbZ$ such that $x^2 - aby^2 = -k$. This implies
that $-k \equiv x^2\,  (\mathrm{mod}\,p)$ for any prime $p$ dividing $a$.
For any $k\in \{1,\ldots, n\}$ choose a prime $p_k$ that violates the
condition $-k \equiv x^2\,  (\mathrm{mod}\,p_k)$ for any
$x$ ($-k$ is not a quadratic residue).
By \ref{lem_integers} below there exist infinitely many such primes,
so we may also assume that $p_k$ are distinct and greater than $b$. Then define $a = \prod_k p_k$.
By construction, the integer $ab$ is not a square, hence $\bbL_{ab}$ does not
represent zero, and since there are infinitely many choices of $p_k$ for every $k$,
there are also infinitely many choices of $a$ as above.
\endproof

\hfill

\lemma\label{lem_integers}
For any integer $k > 0$ there exist infinitely many primes $p$ such that
$-k$ is not a quadratic residue modulo $p$.

\hfill

\proof We have $k = 2^{a_0}\prod_i q_i^{a_i}$ where $q_i$ are distinct odd primes, $a_i > 0$. Then
\begin{equation}\label{eqn_Legendre}
	\left(\frac{-k}{p}\right) = \left(\frac{-1}{p}\right) \left(\frac{2}{p}\right)^{a_0} \prod_{i} \left(\frac{q_i}{p}\right)^{a_i},
\end{equation}
where $\left(\frac{a}{b}\right)$ denotes the Legendre symbol.
Imposing the condition $p\equiv -1\, (\mathrm{mod}\,8)$ will guarantee that
$$
\left(\frac{-1}{p}\right) = -1, \quad \left(\frac{2}{p}\right) = 1,
$$
implying that the product of the first two factors in (\ref{eqn_Legendre}) is equal to $-1$
for any $a_0$.
Moreover, assuming that $p \neq q_i$ for all $i$, the congruence $p\equiv -1\, (\mathrm{mod}\,8)$ implies that we have
$$
\left(\frac{q_i}{p}\right) = (-1)^{\frac{q_i-1}{2}}\left(\frac{p}{q_i}\right),
$$
by the quadratic reciprocity. Choosing the residues of $p$ modulo $q_i$ 
we can guarantee that $\left(\frac{q_i}{p}\right) = 1$ for all $q_i$.
Then (\ref{eqn_Legendre}) will imply that $\left(\frac{-k}{p}\right) = -1$.
By Dirichlet's theorem on arithmetic progressions there exist infinitely many primes $p$ that
satisfy the above congruences.
\endproof

\hfill

\remark If one needs to construct a binary Lorentzian anisotropic lattice $\bbM$ that
satisfies the condition $\mu(\bbM) < -n$ for a given $n$, then the lattices
from \ref{prop_bin_lattices} may be not the optimal choice from the computational
point of view, in the sense that their discriminants are unnecessarily huge.
In fact, a more practical solution may be to consider $\bbM = \bbL_{a^2 - 1}$
for $a > 1$. One can check by a direct computation, using the continued
fraction expansion of $\sqrt{a^2 - 1}$, that $\mu(\bbL_{a^2 - 1}) = 2 - 2a$, and then one can
take any $a > 1 + n /2$.

\subsection{Divisors in the period domain}\label{sec_div}

Given a lattice of rank two $\bbL\subset \Lambda$ that contains the given
element $h\in \Lambda$, the manifold $\DD_{\bbL}$ is a divisor in the period
domain $\DD_h$. We will construct a family of such divisors using the lattices
from Section \ref{sec_binlat}. In what follows, we will assume that $r = \rk(\Lambda)\ge 6$.
Recall that $h$ is a primitive vector that defines the polarization and $d = q(h) > 0$.

Let $\Lambda_h\subset \Lambda$ be the orthogonal complement of $h$. Then $\Lambda_h$ is
a lattice of signature $(2, r - 3)$, therefore it
represents zero (\cite{_Meyers_}). Let $e\in \Lambda_h$
be a primitive isotropic vector, and let $(e, \Lambda_h) = m\bbZ$. Then $\Lambda_h$
and $\tilde{e} = e / m$ span an integral lattice $\widetilde{\Lambda}_h$ in the vector
space $\Lambda_h\otimes \bbQ$. We have an inclusion $\Lambda_h \subset \widetilde{\Lambda}_h$
of index $m$. Now, $(\tilde{e}, \widetilde{\Lambda}_h) = \bbZ$, hence we can find an
isotropic element $\tilde{f}\in \widetilde{\Lambda}_h$ with $(\tilde{e}, \tilde{f}) = 1$.

Given an integer $a$, let $u_a = \tilde{e} - da\tilde{f}$. Then $u_a\in \widetilde{\Lambda}_h$
is a primitive vector with $q(u_a) = -2da$, and $\bbZ u_a\oplus \bbZ h$ is a primitive sublattice
in $\widetilde{\Lambda}_h \oplus \bbZ h$. By construction, $\bbZ u_a\oplus \bbZ h \simeq \bbL_{2a}(d)$.

Recall from Section \ref{sec_hk} that we have a collection of MBM classes in
$\Lambda$ and that there exists $N > 0$ such that for any $v\in\MBM$
we have $-N \le q(v) < 0$.

\hfill

\proposition\label{prop_mj}
There exists an infinite sequence of primitive binary aniso\-tropic sublattices
$\bbM_j\subset \Lambda$, $j\in \bbN$,
with $h\in \bbM_j$ and $\bbM_j \cap \Lambda_h = \bbZ v_j$, where $v_j$ is primitive
with $q(v_j)\to -\infty$ for $j\to +\infty$, and such that $\bbM_j \cap \MBM = \emptyset$ for all $j\in \bbN$.

\hfill

\proof Consider the lattices $\bbZ u_a\oplus \bbZ h \simeq \bbL_{2a}(d)$ constructed above.
Denote by $T$ the index of $\Lambda_h \oplus \bbZ h$ in $\Lambda$. Let $a_j\to +\infty$ be a sequence of integers obtained from \ref{prop_bin_lattices} with the constant $n = NT^2$,
where $N$ is the bound for the squares of the MBM classes, as recalled above.

We define $v_j$ to be a generator of $\bbZ u_{a_j} \cap \Lambda_h$.
Since $[\widetilde{\Lambda}_h: \Lambda_h] = m$, we have $v_j = m_j u_{a_j}$, where $m_j$ is a divisor of $m$, so $q(v_j)\to -\infty$. By construction $\bbZ v_j\oplus \bbZ h$ is
a sublattice of $\bbL_{2a_j}(d)$, so it does not represent the integers $-dNT^2,\ldots,-1,0$.
Let $\bbM_j$ be the primitive sublattice of $\Lambda$ spanned by $\bbZ v_j\oplus \bbZ h$.
Then the index of the latter in $\bbM_j$ is not bigger than $T$. Therefore
$\mu(\bbM_j) < -dN$, because if $x\in \bbM_j$, then $Tx\in \bbZ v_j\oplus \bbZ h$,
and $q(Tx) = T^2q(x)$. We conclude that the lattices $\bbM_j$ do not contain MBM classes.
They are anisotropic, because they are commensurable with $\bbL_{2a_j}(d)$.
Also, by construction, $\bbM_j\cap \Lambda_h = \bbZ v_j$, with $q(v_j)\to -\infty$.
\endproof

\hfill

We now use the lattices constructed above to define a family of divisors in the period domain
and in the moduli space.

\hfill

\definition\label{def_divisors}
Let $\bbM_j$ be the lattices from \ref{prop_mj} and $\DD_{\bbM_j}\subset \DD_h$ the corresponding
divisors. Note that $\DD_{\bbM_j} = \DD_h \cap \bbP(v_j^\perp\otimes \bbC)$, where $v_j$ is the
sequence of vectors appearing in the above proposition. We define $\YY_j$ to be the images of $\DD_{\bbM_j}$
in $\MM_h$ under the quotient map $\DD_h \to \MM_h$.

\subsection{Strong density of a family of divisors}\label{sec_density}

We next discuss the density properties of the divisors constructed above.

\hfill

\definition\label{def_density}
We call a family of divisors $\{\YY_j\}_{j\in \bbN}$ in a complex analytic space $\MM$
strongly dense if the following condition holds: for any
non-constant holomorphic map $f\colon \Delta\to \MM$, where $\Delta$
is the unit disc in $\bbC$, the subset $f^{-1}(\cup_j\YY_j)$ is dense in $\Delta$.

\hfill

\remark
Strong density clearly implies that $\cup_j \YY_j$ is dense in $\MM$
in the analytic topology, because through any point $p\in \MM$ there
pass arbitrarily small holomorphic discs, and any of those discs has
to intersect at least one of the divisors $\YY_j$. The converse is not 
true. For example, consider the linear projection $\pi\colon \bbC P^2 \dashrightarrow \bbC P^1$
from a fixed point in $\bbC P^2$. For a countable and dense set of points
in $\bbC P^1$, their preimages under $\pi$ form a family of lines that
is dense but not strongly dense in $\bbC P^2$.

\hfill

Our proof of strong density of the divisors constructed in \ref{def_divisors} is similar to Oguiso's arguments from \cite{_Oguiso_}: to find an intersection with $\{\YY_j\}$ in a neighbourhood of a point $f(0)$, we first draw a hyperplane through this point, in such a way that it is transversal to $f(\Delta)$, and then approach this hyperplane by $\{\YY_j\}$. To perform this in our setting
we will need to recall the description of $\DD_h$ as a homogeneous space of 
the orthogonal group, as well as an important step in the proof of the main result of \cite{_AV:Orbits_}.

The quadratic form $q$ induced on the vector space $V = \Lambda_h \otimes \bbR$
has signature $(2, r - 3)$, where $r$ is the second Betti number of our hyperk\"ahler
manifold. Let us fix a subspace $V_0\subset V$ with $\dim(V_0) = 2$ and $q|_{V_0} > 0$.
Choose one of the two isotropic subspaces in $V_0\otimes \bbC$ corresponding to
one of the two points in $\DD_h \cap \bbP(V_0\otimes \bbC)$ and denote it by $\ell_0$ (equivalently, fix an orientation of $V_0$). Then $\DD_h$ is the orbit of $\ell_0$ under the
action of $G = \rmO(V, q) \simeq \rmO(2, r-3)$\footnote{This is best seen when one interprets the period domain as the grassmannian of positive oriented 2-planes in the real second cohomology, indeed to give $l_0$ is the same as to give $V_0$ by taking $V_0$ generated by the real and the imaginary part of $l_0$}. The stabilizer of $\ell_0$ is the compact
group $K$ isomorphic to $\SO(2)\times \rmO(r- 3)$.

Denote by $\Gamma\subset G$ the image of $\Mon_h(M)$ under the natural action on the orthogonal complement of $h$. By the above description we see that
$\DD_h\simeq G / K$ and $\MM_h \simeq \Gamma \backslash G / K$.
Choose a point $v_0\in V_0^\perp$ with $q(v_0) = -1$ and consider its stabilizer
$P\subset G$. Then $P\simeq \rmO(2, r-4)$ and the $P$-orbit of $\ell_0$ in $\DD_h$
is a divisor $\DD_0$ isomorphic to $P/ (P\cap K)$ as a homogeneous space.
There exists a sequence of elements $\gamma_j\in G$, $j\in \bbN$, such that
$\gamma_j v_0 = \lambda_j v_j$, where $v_j$ are the vectors appearing in \ref{def_divisors},
and $\lambda_j\in \bbR$. As for the divisors appearing in that
definition, we have $\DD_{\bbM_j} = \gamma_j \DD_0$. In other words, $\DD_{\bbM_j}$
and $\YY_j$ are the images of the left coset $\gamma_j P\subset G$ under the projection
$G\to G/K$, respectively $G\to \Gamma\backslash G/ K$.

The group $\Gamma$ is countable and we may enlarge the set of vectors $v_j$ appearing
in \ref{def_divisors} adding the $\Gamma$-orbit for each of those vectors.
Then we may assume that the set $\{v_j\}_{j \in \bbN}$ is stable under the action of $\Gamma$.
The condition on the squares of the vectors $v_j$ appearing in \ref{prop_mj}
guarantees that the set $\{v_j\}_{j \in \bbN}$ consists of infinitely many $\Gamma$-orbits,
and therefore we also have infinitely many $\Gamma$-orbits in the set $\{\gamma_j\}_{j\in \bbN}$
of the elements appearing above.

\hfill

\theorem \cite{_AV:Orbits_}\label{thm_AV} In the above setting the union of the
left cosets $\bigcup_{j\in \bbN} \gamma_j P$ is dense in $G$ (and $\gamma_j P$ are dense in the space of left cosets).

\hfill

\proof
Although not claimed explicitly, this is the main step of the proof of \cite[Theorem 1.7]{_AV:Orbits_} applied to our situation.
\endproof

\hfill

In particular there is a sequence of $\DD_{\bbM_j}$ converging to any translate $g_0D_0$. It remains to show that we can choose $g_0$ in such a way that $g_0D_0$ is transverse to a given image of a disc.

For this we analyze the action of $G$ on $\bbP(T^*\DD_h)$, where $T^*\DD_h$ is the holomorphic
cotangent bundle. Let us choose a base point $\tilde{\ell}_0\in \bbP(T^*\DD_h)$ lying above
$\ell_0\in \DD_h$. To do this, observe that $T_{\ell_0}\DD_h$ is naturally
isomorphic to $V_0^\perp\otimes \bbC$. The orthogonal complement of $v_0$ inside $T_{\ell_0}\DD_h$
is a hyperplane, and we let $\tilde{\ell}_0$ to be the corresponding point in the
projectivisation of $T^*_{\ell_0}\DD_h$. Denote by $\widetilde{\DD}_h$ the $G$-orbit
of $\tilde{\ell}_0$.

The $P$-orbit of $\tilde{\ell}_0$, denoted by $\widetilde{\DD}_0$, is isomorphic to $P/(P\cap K)\simeq \DD_0$,
and this orbit is the canonical lifting of the divisor $\DD_0$ to the projectivised
cotangent bundle of $\DD_h$, where each point of the divisor is lifted to the
hyperplane (considered as a point in the dual space) tangent to the divisor.

Let us now consider the intersection $\widetilde{\DD}_h\cap \bbP(T^*_{\ell_0}\DD_h)$.
This intersection is the $K$-orbit of $\tilde{\ell}_0$. Note that 
$K \simeq \SO(2)\times \rmO(r- 3)$, and the first factor acts trivially
on $V_0^\perp\otimes \bbC\simeq T_{\ell_0}\DD_h$, while for the second factor
the latter space is the complexification of the standard representation.
In particular, it follows that $T^*_{\ell_0}\DD_h$ is irreducible as a complex
representation of $K$.

\hfill

\proposition\label{prop_density}
In the above setting, the collection of divisors $\{\gamma_j\DD_0\}_{j\in \bbN}$
is strongly dense in $\DD_h$, and $\{\YY_j\}_{j\in \bbN}$ is strongly dense in $\MM_h$.

\hfill

\proof Let $f\colon \Delta\to \DD_h$ be a non-constant holomorphic map.
Since the choice of the base point $\ell_0\in \DD_h$ was arbitrary, we may
assume that $\ell_0 = f(0)$. Consider the tangent vector $f'(0)\in T_{\ell_0}\DD_h$.
As we have observed above, the vector space $T^*_{\ell_0}\DD_h$ does not contain any
$K$-invariant complex subspaces. It follows that there exists $g_0\in K$ such that
the hyperplane in $T_{\ell_0}\DD_h$ defined by $g_0\tilde{\ell}_0$ does not contain
$f'(0)$. In this case the divisor $g_0\DD_0$ passes through the point $f(0)$ but
does not contain the image of $f$. By the density result of \ref{thm_AV}, the divisor
$g_0\DD_0$ is the limit of a sequence of divisors $\gamma_{j_\alpha}\DD_0$, $\alpha\to +\infty$.
Note that all of those divisors are intersections of $\DD_h$ and some
hyperplanes in $\bbP(V\otimes \bbC)$. Our claim about the divisors $\{\gamma_j\DD_0\}_{j\in \bbN}$
now follows from \ref{lem_dens} below. The claim about the divisors $\YY_j$
follows as well, because $\YY_j$ are the images of $\gamma_j\DD_0$ under
the quotient map $\DD_h \to \DD_h/ \Gamma$.
\endproof

\hfill

\lemma\label{lem_dens}
Let $f\colon \Delta\to \bbC P^n$ be a non-constant holomorphic map
and $H\subset \bbC P^n$ be a hyperplane such that
$f(0)\in H$ but $f(\Delta)\nsubset H$. Then for any sequence of
hyperplanes $H_i$, $i\in \bbN$ that converge to $H$
(as points of the dual projective space) when $i\to +\infty$, there
exists $i_0\in \bbN$ such that $f(\Delta)\cap H_{i_0} \neq \emptyset$.

\hfill

\proof
We may pass to a local chart containing $f(0)$ and 
assume that $H_i$ are defined by degree one polynomials $g_i$ that
converge to a degree one polynomial $g$ that defines $H$. Then $\varphi = g\circ f$
is a non-constant holomorphic function on $\Delta$ with $\varphi(0) = 0$,
and $\varphi_i = g_i\circ f$ are holomorphic functions converging to
$\varphi$ uniformly on compact neighborhoods of $0\in \Delta$.
The claim now follows from Rouché's theorem.
\endproof

\section{Proofs}\label{sec_proofs}

\subsection{Proof of \ref{thm_main}}\label{sec_proofthm}

\begin{enumerate}
	\item The family of subvarieties $\XX_i$ is defined in Section \ref{sec_nikvin},
	see \ref{def_NV}. Using the notations from the proof of \ref{prop_NV}, we see
	that $\NS(M)\in \NN\VV^+$. Then the orthogonal complement of $h$ in $\NS(M)$
	is one of the lattices $\bbM_i$ constructed in that proposition.
	The period of $M$ is contained in the special subvariety $\DD_{\bbM_i}$,
	and for a very general point $[M']\in \DD_{\bbM_i}$ we have $\NS(M')\simeq \bbM_i$.
	This completes the proof of the first part.
	\item The family of divisors $\{\YY_j\}_{j\in \bbN}$ is given by \ref{def_divisors}.
	Their strong density is \ref{prop_density}.
	\item The hyperk\"ahler manifolds parametrized by the divisors $\YY_j$ are projective,
	and there exist universal families $\pi_j\colon \UU_j\to \YY^\circ_j$ over some non-empty
	Zariski open subsets $\YY_j^\circ\subset \YY_j$ (possibly after passing to a finite covering
	of those subsets). For a very general point $t\in \YY_j^\circ$, the fibre $M = \pi_j^{-1}(t)$
	has $\NS(M)$ isomorphic to one of the rank two lattices constructed in \ref{prop_bin_lattices}.
	This lattice, by construction, contains no MBM classes, is indefinite and does not
	represent zero. Hence its orthogonal group contains an infinite cyclic subgroup,
	and the same is true about
	the automorphism group of $M$, see Section \ref{sec_hk}. So $\Aut(M)$ contains an element
	of infinite order when $[M]\in \YY_j^\circ$ is a very general point. Using the standard
	Hilbert scheme argument (using the fact that the family $\UU_j$ is projective),
	we may shrink $\YY_j^\circ$ and guarantee that the same holds for any $[M]\in \YY_j^\circ$.
\end{enumerate}

\subsection{Proof of \ref{cor_main}}\label{sec_proofcor}

After possibly passing to the universal cover of $B$ we may
consider the period map $\psi\colon B\to \MM_h$. Let $\bbM$ be
the lattice that represents the isomorphism class of $\NS(M_t)$
for a very general $t\in B$. If the lattice $\bbM$ is not
reflective, then $\Bim(M_t)$ is infinite for a very general $t\in B$,
completing the proof. If $\bbM$ is reflective, then we consider three cases.

First assume that $\rk(\bbM)\ge 3$.  Since $\bbM$ is reflective, $\mathrm{Im}(\psi)\subset \XX_{i_0}$
for one of the subvarieties $\XX_{i_0}$ from the first part of \ref{thm_main}.
We may also choose $\XX_{i_0}$ to be the of the smallest possible dimension among $\XX_i$ containing $\mathrm{Im}(\psi)$.
Then, since our family is not isotrivial, we can find a dense subset $B''\subset B$ such
that $\rk \NS(M_t) > \rk \bbM$ for $t\in B''$ (see e.g. \cite{_Oguiso_}). Let 
$$J = \{j\st \dim(\XX_j) < \dim(\XX_{i_0})\}.$$
The image $\psi(B'')\subset \XX_{i_0}$ can not be contained in the union $\cup_{j\in J} \XX_j$, indeed otherwise it must be contained in one of $\XX_j, j\in J$, and this is not possible by the choice of $i_0$.
Setting $B' = \{t\in B''\st \psi(t)\notin \cup_{j\in J} \XX_j\}$,
we get the desired dense subset in $B$.

The second case is when $\rk(\bbM)= 2$. In this case we 
can find a dense subset $B''\subset B$ such that $\rk \NS(M_t) > 2$ for $t\in B''$.
As above, $\psi(B'')$ can not be contained in the Nikulin--Vinberg locus $\cup_{i} \XX_i$,
and we conclude by defining $B' = \{t\in B''\st \psi(t)\notin \cup_{i} \XX_i\}$.

Finally, if $\rk(\bbM)= 1$, by the strong density of the divisors $\YY_j$
from \ref{thm_main}, the preimage $B'' = \psi^{-1}(\cup_j \YY_j)$ is dense in $B$.
For any point $t\in B''$, either $\NS(M_t)$ is of rank two and non-reflective (by the
definition of $\YY_j$), or $\rk \NS(M_t)\ge 3$. Then, as above, we conclude 
by defining $B' = \{t\in B''\st \psi(t)\notin \cup_{i} \XX_i\}$.

\hfill



{\small
\noindent {\sc Ekaterina Amerik\\
{\tt  Ekaterina.Amerik@gmail.com}\\
{\sc  Universit\'e Paris-Saclay,\\
Laboratoire de Math\'ematiques,\\
Campus Scientifique d'Orsay, B\^atiment 307, 91405 Orsay, France}\\
{\sc Laboratory of Algebraic Geometry,\\
National Research University HSE,\\
Department of Mathematics, 6 Usacheva Str. Moscow, Russia}

\hfill

\noindent {\sc Andrey Soldatenkov\\
{\tt aosoldatenkov@gmail.com}\\
{\sc Universidade Estadual de Campinas\\
Departamento de Matem\'atica - IMECC\\
Rua S\'ergio Buarque de Holanda, 651\\
13083-859, Campinas - SP, Brasil}}

\hfill

\noindent {\sc Misha Verbitsky\\
{\sc Instituto Nacional de Matem\'atica Pura e
              Aplicada (IMPA) \\ Estrada Dona Castorina, 110\\
Jardim Bot\^anico, CEP 22460-320\\
Rio de Janeiro, RJ - Brasil\\
\tt  verbit@impa.br }
}}}

\end{document}